\documentclass[11pt]{article}
\usepackage{amsmath}
\usepackage{amssymb}
\usepackage{amsthm}
\usepackage{url}

\newcommand\mLP{\\[\medskipamount]}
\newcommand\bLP{\\[\bigskipamount]}
\newcommand\CC{\mathbb{C}}
\newcommand\RR{\mathbb{R}}
\newcommand\ZZ{\mathbb{Z}}
\newcommand\al\alpha
\newcommand\be\beta
\newcommand\ga\gamma
\newcommand\de\delta
\newcommand\ep\varepsilon
\newcommand\la\lambda
\newcommand\De{\Delta}
\newcommand\La{\Lambda}
\newcommand\half{\frac12}
\newcommand\thalf{\tfrac12}
\newcommand\iy\infty
\newcommand\lan{\langle}
\newcommand\ran{\rangle}
\newcommand\Znonneg{\ZZ_{\ge0}}
\newcommand{\qhyp}[5]{\,\mbox{}_{#1}\phi_{#2}\left(
  \genfrac{}{}{0pt}{}{#3}{#4};#5\right)}
\newcommand\LHS{left-hand side}
\newcommand\RHS{right-hand side}
\numberwithin{equation}{section}
\newtheorem{theorem}{Theorem}[section]
\newtheorem{proposition}[theorem]{Proposition}
\newtheorem{Remark}[theorem]{Remark}
\newenvironment{remark}{\begin{Remark}\rm}{\end{Remark}}
\newcommand\Proof{\smallskip\noindent{\bf Proof}\quad}

\begin{document}
\title{The structure relation for Askey-Wilson polynomials}
\author{Tom H. Koornwinder\bLP
{\small  Korteweg-de Vries Institute, University of Amsterdam,}\\
{\small Plantage Muidergracht 24, 1018 TV Amsterdam, The Netherlands}\\
{\small {\tt thk@science.uva.nl}}
}
\date{\small{This paper is dedicated to Nico Temme
on the occasion of his 65th birthday.}}
\maketitle
\begin{abstract}
An explicit structure relation for Askey-Wilson polynomials is given.
This involves a divided $q$-difference
operator which is skew symmetric with respect to the
Askey-Wilson inner product and which sends polynomials of degree $n$ to
polynomials of degree $n+1$. By specialization of parameters and by taking
limits, similar structure relations, as well as lowering and raising
relations, can be obtained for other families
in the $q$-Askey scheme and the Askey scheme. This is explicitly discussed
for Jacobi polynomials, continuous $q$-Jacobi polynomials,
continuous $q$-ultraspherical polynomials,
and for big $q$-Jacobi polynomials.
An already known structure relation for this last family
can be obtained from the new
structure relation by using the three-term recurence relation and
the second order $q$-difference formula.
The results are also put in the framework of a more general theory.
Their relationship with earlier work by Zhedanov and Bangerezako is discussed.
There is also a connection with the string equation in discrete matrix models
and with the Sklyanin algebra. 
\end{abstract}
% until 78
\section{Introduction}
One of the many ways to characterize a family $\{p_n(x)\}$
of classical orthogonal polynomials
(Jacobi, Laguerre and Hermite polynomials) is by
a {\em structure relation}
\begin{equation}
\pi(x)p_n'(x)=a_np_{n+1}(x)+b_np_n(x)+c_np_{n-1}(x),
\label{eq:20}
\end{equation}
where $\pi(x)$ is a fixed polynomial (necessarily of degree $\le 2$).
In view of the three-term recurrence relation this is equivalent to
a characterization by a {\em lowering relation}
\begin{equation}
\pi(x)\,p_n'(x)=(\al_n x+\be_n)p_n(x)+\ga_n p_{n-1}(x)
\label{eq:21}
\end{equation}
or to a characterization by a {\em raising relation}
\begin{equation}
\pi(x)\,p_n'(x)=(\tilde\al_n x+\tilde\be_n)p_n(x)+\tilde\ga_n p_{n+1}(x).
\label{eq:22}
\end{equation}
The lowering relation \eqref{eq:21} was given in
\cite[10.7(4)]{01} (in \cite[10.8(15)]{01} explicitly for Jacobi polynomials),
where it was attributed to Tricomi (1948). The characterization of
classical orthogonal polynomials by their property \eqref{eq:21}
(and thus equivalently by their property \eqref{eq:20})
was first given by Al-Salam \& Chihara \cite{02}. See also
\cite{07} for a characterization by \eqref{eq:20}.
Note that the lowering and raising relations \eqref{eq:21}, \eqref{eq:22}
are different from the more familiar {\em shift operator relations},
where not only the degree is lowered or raised, but also the parameters
are shifted. These formulas are well known for all orthogonal
polynomials in the \mbox{($q$-)Askey} scheme, see for instance
\cite{04}.

Orthogonal polynomials satisfying the more general structure relation
\begin{equation}
\label{eq:06}
\pi(x)\,p_n'(x)=\sum_{j=n-s}^{n+t}a_{n,j} p_j(x)\quad
\mbox{($\pi(x)$ a polynomial; $s,t$ independent of $n$).}
\end{equation}
are called {\em semi-classical}, see Maroni \cite{03}.
According to \cite{02} the question to characterize orthogonal polynomials
satisfying \eqref{eq:06} was first posed by R.~Askey.

Garc{\'\i}a, Marcell\'an \& Salto \cite{08} characterized discrete
classical orthogonal polynomials in the Hahn class (Hahn, Krawtchouk, Meixner
and Charlier polynomials) by a structure relation similar to \eqref{eq:20},
with the derivative replaced by the difference operator
$(\De f)(x):=f(x+1)-f(x)$.
Next Medem, \'Alvarez-Nodarse \& Marcell\'an \cite{09}
(see also \cite{10})
characterized the orthogonal polynomials in the $q$-Hahn class by
a structure relation obtained from \eqref{eq:20} by replacing the derivative
by the $q$-derivative $D_q$, where
\begin{equation}
(D_qf)(x)=D_{q,x}\bigl(f(x)\bigr):=\frac{f(x)-f(qx)}{(1-q)x}.
\label{eq:39}
\end{equation}
Here a family of orthogonal polynomials $p_n(x)$ is in the {\em $q$-Hahn class}
if the polynomials $(D_qp_n)(x)$ are again orthogonal.

Variants of lowering and raising relations \eqref{eq:21}, \eqref{eq:22}
are scattered over the literature. See a brief survey in \cite{11}.
Note in particular the lowering and raising relations for $A_n$ type
Macdonald polynomials given by Kirillov \& Noumi \cite{12}.
The $A_1$ case yields lowering and raising relations for
continuous $q$-ultraspherical polynomials.
A lowering relation for continuous $q$-Jacobi polynomials was given by
Ismail \cite[Theorem 15.5.2]{24}.

For Askey-Wilson polynomials, a structure relation was given in a very
implicit way by Zhedanov \cite{16} (see my discussion in Remark \ref{th:78}),
while a lowering and raising relation was given by
Bangerezako \cite{20} (see my discussion after \eqref{eq:77}).
The main result of the present paper
gives a structure relation for Askey-Wilson polynomials in the form
\begin{equation}
Lp_n=a_np_{n+1}+c_np_{n-1},
\label{eq:23}
\end{equation}
where $L$ is a linear operator which is  skew symmetric with respect to the
inner product $\lan .,.\ran$
for which the Askey-Wilson polynomials are orthogonal:
\begin{equation}
\lan Lf,g\ran=-\lan f,Lg\ran.
\label{eq:24}
\end{equation}
This operator is explicitly given by
\begin{multline}
(Lf)[z]:=
\Bigl((1-az)(1-bz)(1-cz)(1-dz)\,z^{-2}\,f[qz]\\
-(1-a/z)(1-b/z)(1-c/z)(1-d/z)\,z^2\,f[q^{-1}z]\Bigr)\,(z-z^{-1})^{-1}.
\label{eq:19}
\end{multline}
It sends symmetric Laurent polynomials of degree $n$ to
symmetric Laurent polynomials
of degree $n+1$. By specialization and limit transition a structure relation
of the form \eqref{eq:23} with $L$ satisfying \eqref{eq:24} can be obtained
for all families of orthogonal polynomials in the $q$-Askey scheme and
the Askey scheme (a list of these families is given in \cite{04}).
For instance, for Jacobi polynomials $P_n^{(\al,\be)}(x)$ we get
for the operator $L$:
\begin{align}
&(Lf)(x):=(1-x^2)f'(x)-\thalf\bigl(\al-\be+(\al+\be+2)x\bigr)f(x)
\label{eq:29}
\\
&\quad=\,(1-x)^{-\half\al+\half}(1+x)^{-\half\be+\half}\,
\frac d{dx}\Bigl((1-x)^{\half\al+\half}(1+x)^{\half\be+\half}\,f(x)\Bigr)
\label{eq:38}
\\
&\quad=((DX-XD)f)(x),
\label{eq:68}
\end{align}
where $X$ is multiplication by $x$ and $D$ is a second order differential
operator having the Jacobi polynomials as eigenfunctions
(see \eqref{eq:72}).
It will turn out that the form \eqref{eq:38} of $L$, involving something
close to the square root of the weight function, can also be realized
higher up in the \mbox{($q$-)Askey} scheme, notably in the Askey-Wilson case
\eqref{eq:19}.
It will also turn out that the form \eqref{eq:68} of $L$, i.e., as
a commutator of $X$ and a second order differential or
\mbox{($q$-)difference} operator $D$ having the orthogonal polynomials as
eigenfunctions, persists in the \mbox{($q$-)Askey} scheme.
In fact, there is an essentially one-to-one relationship between
operators $L$ and $D$.

For those families where a structure relation had been given earlier, one
can relate that formula to \eqref{eq:23} by use of the three-term recurrence
relation, and sometimes also of the second-order
\mbox{($q$-)difference} equation.

The general theory of the structure relation of the form \eqref{eq:23}
with skew symmetric $L$ will be given in \S\ref{sec:43}. This theory is
easy and elegant. The coefficients in the resulting structure relation
and in the lowering and raising relations
are very close to the coefficients in the
three-term recurrence relation. In \S\ref{sec:43} I will also discuss the
relationship with bispectral problems (Gr\"unbaum and Haine),
Zhedanov's algebra, and the string equation in the context of discrete
matrix models.
The case of Jacobi polynomials
will be discussed in \S\ref{sec:44}. The main result, the Askey-Wilson
case, is the topic of \S\ref{sec:36}. Here also a connection with the
Sklyanin algebra will be made.
In \S\ref{sec:45} this is specialized to the case of continuous
$q$-Jacobi polynomials and we show that it has the results for
Jacobi polynomials as a limit case.
A further specialization to
continuous $q$-ultraspherical polynomials is given in
\S\ref{sec:47}, and the resulting structure relation is related to
another one obtained from results in \cite{11}.
Finally, in \S\ref{sec:46}, we take the limit of the Askey-Wilson case
to the case of big $q$-Jacobi polynomials, and we relate the resulting
structure relation to the one in \cite{09}.
\bLP
{\bf Acknowledgement}\quad
I thank
G. Bangerezako,
L. Haine,
H. Rosengren,
P.~van Moerbeke
and the referees for helpful comments.
\bLP
{\bf Conventions}\\
Throughout assume that $0<q<1$. For \mbox{($q$-)Pochhammer} symbols and
\mbox{($q$-)hypergeometric series}
use the notation of Gasper \& Rahman \cite{06}.
Symmetric Laurent polynomials
$f[z]=\sum_{k=-n}^n c_k z^k$ (where $c_k=c_{-k}$) are related to ordinary
polynomials $f(x)$ in $x=\thalf(z+z^{-1})$ by
$f(\thalf(z+z^{-1}))=f[z]$.
\section{The general form of the structure relation}
\label{sec:43}
Suppose we have a family of orthogonal polynomials $p_n(x)$ with
respect to an
orthogonality measure $\mu$ on $\RR$:
\begin{equation}
p_n(x)=k_n x^n+\cdots\,,\qquad
\lan p_n,p_m\ran=h_n\de_{n,m},
\end{equation}
where
\begin{equation}
\lan f,g\ran:=\int_\RR f(x)\,g(x)\,d\mu(x).
\label{eq:25}
\end{equation}
Write the three-term recurrence relation as
\begin{equation}
xp_n(x)=A_n p_{n+1}(x)+B_n p_n(x)+C_np_{n-1}(x).
\label{eq:58}
\end{equation}
Then
\begin{equation}
A_n=\frac{k_n}{k_{n+1}}\,,\qquad
C_n=\frac{k_{n-1}}{k_n}\,\frac{h_n}{h_{n-1}}=A_{n-1}\,\frac{h_n}{h_{n-1}}\,.
\end{equation}
The proof of the following proposition is straightforward.
\begin{proposition}
\label{th:27}
Let $L$ be a linear operator acting on the space $\RR[x]$
of polynomials in one
variable with real coefficients such that $L$ is skew symmetric with respect
to the inner product \eqref{eq:25} (i.e., \eqref{eq:24} holds) and such that
\begin{equation}
L(x^n)=\ga_n x^{n+1}+\mbox{\rm terms of lower degree},
\label{eq:57}
\end{equation}
where $\ga_n\ne0$. Then the following structure relation and
lowering and raising relations hold:
\begin{align}
(Lp_n)(x)=&\,\ga_nA_np_{n+1}(x)-\ga_{n-1}C_np_{n-1}(x),
\label{eq:28}\mLP
\noalign{\allowbreak}
-\ga_n(x-B_n)p_n(x)+(Lp_n)(x)=&\,-(\ga_n+\ga_{n-1})C_np_{n-1}(x),
\label{eq:31}\\
\noalign{\allowbreak}
\ga_{n-1}(x-B_n)p_n(x)+(Lp_n)(x)=&\,(\ga_n+\ga_{n-1})A_np_{n+1}(x).
\label{eq:32}
\end{align}
\end{proposition}

Skew symmetric operators $L$ as above can be produced from symmetric operators
$D$ which have the $p_n$ as eigenfunctions. First note that
\begin{equation}
(Xf)(x):=xf(x)
\label{eq:70}
\end{equation}
defines a symmetric operator $X$ on $\RR[x]$ with respect to the
inner product \eqref{eq:25}, i.e., $\lan Xf,g\ran=\lan f,Xg\ran$.
Now the following proposition can be
shown immediately.
\begin{proposition}
\label{th:62}
Let $D$ be a linear operator acting on $\RR[x]$ which is symmetric
with respect to the inner product \eqref{eq:25}, i.e.,
$\lan Df,g\ran=\lan f,Dg\ran$, and which satisfies
\begin{equation}
D(x^n)=\la_n x^n+\mbox{\rm terms of lower degree},
\label{eq:63}
\end{equation}
where $\la_n\ne \la_{n-1}$, so $p_n$ is an eigenfunction of $D$:
\begin{equation}
Dp_n=\la_np_n.
\label{eq:64}
\end{equation}
Then the commutator 
\begin{equation}
L:=[D,X]=DX-XD
\label{eq:65}
\end{equation}
is skew symmetric with respect to the inner product \eqref{eq:25}
and satisfies \eqref{eq:57} with
\begin{equation}
\ga_n=\la_{n+1}-\la_n\ne0.
\label{eq:66}
\end{equation}
So $L$ also satisfies the structure relation \eqref{eq:28}.
\end{proposition}

In fact, we can reverse Proposition \ref{th:62}: From skew symmetric
$L$ as in Proposition \ref{th:27} we can produce symmetric $D$
as in Proposition \ref{th:62}, and $D$ is uniquely determined by $L$
up to a term which is constant times identity.
\begin{proposition}
Let $L$ be as in Proposition \ref{th:27}. Define a linear operator $D$ on
$\RR[x]$ by its action on monomials:
\begin{equation}
D(1)=0,\quad
D(x^n)=\sum_{k=0}^{n-1} X^kL(x^{n-k-1}).
\label{eq:67}
\end{equation}
Then $D$ satisfies the properties of Proposition \ref{th:62} with $\la_0=0$.
Any other operator $D$ satisfying \eqref{eq:65} and having 1 as
eigenfunction differs from $D$
given by \eqref{eq:67} by a constant times identity.
\end{proposition}
\Proof
Formula \eqref{eq:65} acting on $x^n$ follows directly from
\eqref{eq:67}. From \eqref{eq:57} together with $DX=L+XD$ acting on $x^n$
we see by induction that \eqref{eq:63} holds with $\la_n$ satisfying
\eqref{eq:66}. In order to prove that $D$ is symmetric, observe that,
for $n>0$,
\begin{align*}
\lan D x^n,x^m\ran-\lan x^n,Dx^m\ran
&=\lan DXx^{n-1},x^m\ran-\lan x^{n-1},XDx^m\ran\\
&=\lan(L+XD)x^{n-1},x^m\ran+\lan x^{n-1},(L-DX)x^m\ran\\
&=\lan Dx^{n-1},x^{m+1}\ran-\lan x^{n-1},Dx^{m+1}\ran.
\end{align*}
Hence
\[
\lan Dx^{n+m},1\ran=\lan Dx^n,x^m\ran-\lan x^n,Dx^m\ran=
-\lan 1,Dx^{n+m}\ran.
\]
So $\lan Dx^{n+m},1\ran=0$ and $\lan Dx^n,x^m\ran=\lan x^n,Dx^m\ran$.

Finally, for the uniqueness, let $D_1$ and $D_2$ satisfy \eqref{eq:67}
and and let them have 1 as eigenfunction.
Then $D_1-D_2$ commutes with $X$. Hence $(D_1-D_2)(x^n)=
X^n(D_1-D_2)(1)=X^n(c\,1)=c\,x^n$.
\qed
\begin{remark}
\label{th:75}
If $L$ satisfies \eqref{eq:65} and $D$ satisfies \eqref{eq:64} then the
lowering and raising relations \eqref{eq:31} and \eqref{eq:32}
are preserved if we add to their left-hand sides a term
$g(x)\,(D-\la_n)p_n(x)$, where $g(x)$ is any function.
\end{remark}
\begin{remark}
Van Moerbeke \cite[\S7]{18} (also jointly with Adler in
\cite[\S5]{17}) gives explicit skew symmetric first order
differential operators $L$ ($Q$ in his notation) satisfying
Proposition 2.1 for the case of Jacobi, Laguerre and Hermite polynomials.
More generally he gives these operators if the weight function is perturbed
by multiplying the weight function with
$\exp(\sum_{j=1}^\iy t_jx^j)$ (only finitely many $t_j$ nonzero).
Then \eqref{eq:57} and \eqref{eq:28} are no longer valid, but the
\RHS\ of \eqref{eq:28} has to be replaced by some linear combination of
orthogonal polynomials $p_k(x)$ with coefficients depending on the $t_j$.
In all these cases $L$ satsifies the so-called {\em string equation}
\begin{equation}
[X,L]=f_0(X)
\label{eq:74}
\end{equation}
with $f_0$ a polynomial given by $f_0(X)=1-X^2$, $X$, $1$
for Jacobi, Laguerre and Hermite polynomials, respectively, and for their
deformations. This is inspired by the matrix models of the physicists,
where the Hermite case occurs, see Witten \cite[\S4c]{19}, in particular
(4.43), (4.53), (4.65), (4.66).
\end{remark}
\begin{remark}
\label{th:78}
The structure relation \eqref{eq:28} can be written in a form which
is symmetric with respect to the variables $n$ and $x$.
Write
$p(n,x)$ instead of $p_n(x)$. Define operators $J$ and $\La$
acting on functions on $\Znonneg$ by
\begin{align}
(J\phi)(n)&:=A_n\phi(n+1)+B_n\phi(n)+C_n\phi(n-1),\\
(\La\phi)(n)&:=\la_n \phi(n).
\end{align}
Now use \eqref{eq:70}, \eqref{eq:65}, \eqref{eq:58} and
\eqref{eq:64} in order to rewrite
\eqref{eq:28} as
\begin{equation}
\label{eq:71}
\bigl([D,X]\,p(n,\,.\,)\bigr)(x)=\bigl([J,\La]\,p(\,.\,,x)\bigr)(n).
\end{equation}
This equation essentially occurs in Duistermaat \& Gr\"unbaum
\cite[(1.7)]{13} and Gr\"unbaum \& Haine \cite[(3)]{15},
where they study the bispectral problem with $D$ a differential operator
and $n$ continuous respectively discrete.

For $p_n$, $D$ and $\la_n$ $q$-dependent, a $q$-analogue of
\eqref{eq:71} has been considered in the literature which involves
$q$-commutators instead of ordinary commutators:
\begin{equation}
\label{eq:73}
\Bigl((q^\half DX-q^{-\half}XD)\,p(n,\,.\,)\Bigr)(x)=
\Bigl((q^\half J\La-q^{-\half}\La J)\,p(\,.\,,x)\Bigr)(n).
\end{equation}
In fact, in Zhedanov \cite{16} formulas (1.4) and (1.8a) may be interpreted
as formulas \eqref{eq:64} and \eqref{eq:58} in the present paper,
respectively. Then formula (1.8b) together with (1.1a) in \cite{16}
can be interpreted as \eqref{eq:73} above.
See also Gr\"unbaum \& Haine \cite{14} on the bispectral problem in
the $q$-case, where some formulas in section 3 may be close to
\eqref{eq:73} above.

Formula (1.1c) in \cite{16} may be interpreted as the
$q$-commutator $q^\half XL-q^{-\half}LX$ being equal to a linear combination
of $X$, $D$ and $L$. A similar formula is observed in
\cite[(2.4)]{14}, where it is called the $q$-string equation as a
$q$-analogue of the string equation \eqref{eq:74}.
\end{remark}
\section{Jacobi polynomials}
\label{sec:44}
{\em Jacobi polynomials}
$P_n(x)=P_n^{(\al,\be)}(x)$
(see \cite[\S10.8]{01}) are orthogonal with respect to the inner product
\[
\lan f,g\ran=\int_{-1}^1 f(x)\,g(x)\,(1-x)^\al(1+x)^\be\,dx
\qquad(\al,\be>-1).
\]
The coefficients in the three-term recurrence relation 
\eqref{eq:58} are
\begin{align}
A_n&=\frac{2(n+1)(n+\al+\be+1)}{(2n+\al+\be+1)(2n+\al+\be+2)}\,,
\\
B_n&=\frac{\be^2-\al^2}{(2n+\al+\be)(2n+\al+\be+2)}\,,
\\
C_n&=\frac{2(n+\al)(n+\be)}{(2n+\al+\be)(2n+\al+\be+1)}\,.
\end{align}
For the operator $L$ given by \eqref{eq:29} we see immediately that
\eqref{eq:57} holds with
\begin{equation}
\ga_n=-\thalf(2n+\al+\be+2).
\end{equation}
and that the skew symmetry \eqref{eq:24} holds.
Hence, by Proposition \ref{th:27}, the structure relation \eqref{eq:28}
is valid. Explicitly it reads as follows:
\begin{multline}
\Bigl((1-x^2)\frac d{dx}-\thalf\bigl(\al-\be+(\al+\be+2)x\bigr)\Bigr)
P_n^{(\al,\be)}(x)\\
=-\,\frac{(n+1)(n+\al+\be+1)}{2n+\al+\be+1}\,
P_{n+1}^{(\al,\be)}(x)
+\frac{(n+\al)(n+\be)}{2n+\al+\be+1}\,
P_{n-1}^{(\al,\be)}(x).
\label{eq:26}
\end{multline}
We can also make explicit formulas 
\eqref{eq:68} and \eqref{eq:66}
with
\begin{equation}
D=\thalf(1-x^2)\frac{d^2}{dx^2}+
\thalf\bigl(\be-\al-(\al+\be+2)x\bigr)\frac d{dx},
\quad
\la_n=-\thalf n(n+\al+\be+1).
\label{eq:72}
\end{equation}

Formula \eqref{eq:26} can be combined
with the three-term recurrence relation
in order to obtain the structure relation of the form \eqref{eq:20}:
\begin{multline}
\label{eq:02}
(1-x^2)\,\frac d{dx} P_n^{(\al,\be)}(x)=
-\,\frac{2n(n+1)(n+\al+\be+1)}{(2n+\al+\be+1)(2n+\al+\be+2)}\,
P_{n+1}^{(\al,\be)}(x)\\
+\frac{2n(n+\al+\be+1)(\al-\be)}{(2n+\al+\be)(2n+\al+\be+2)}\,
P_n^{(\al,\be)}(x)\\
+\frac{2(n+\al)(n+\be)(n+\al+\be+1)}{(2n+\al+\be)(2n+\al+\be+1)}\,
P_{n-1}^{(\al,\be)}(x).
\end{multline}
\section{Askey-Wilson polynomials}
\label{sec:36}
Askey-Wilson polynomials (see \cite{05}, \cite[\S7.5]{06}, \cite[\S3.1]{04})
are defined by
\begin{align}
&p_n[z]=p_n\bigl(\thalf(z+z^{-1})\bigr)
=p_n\bigl(\thalf(z+z^{-1});a,b,c,d\mid q\bigr)
\nonumber\\
&\qquad\qquad:=\frac{(ab,ac,ad;q)_n}{a^n}\,
\qhyp43{q^{-n},q^{n-1}abcd,az,az^{-1}}{ab,ac,ad}{q,q}.
\label{eq:11}
\end{align}
If $a,b,c,d\in\CC$ satisfy
\begin{equation}
a^2,b^2,c^2,d^2,ab,ac,ad,bc,bd,cd\notin\{q^{-k}\mid k=0,1,2,\ldots\}
\end{equation}
then these polynomials satisfy the orthogonality relation
\begin{equation}
\lan f,g\ran:=
\frac1{4\pi i}\oint_C f[z]\,g[z]\,w(z)\,\frac{dz}z,
\label{eq:30}
\end{equation}
\begin{equation}
w(z):=\frac{(z^2,z^{-2};q)_\iy}
{(az,az^{-1},bz,bz^{-1},cz,cz^{-1},dz,dz^{-1};q)_\iy}\,,
\end{equation}
where $C$ is the unit circle traversed in positive direction with suitable
deformations to separate the sequences of poles converging to zero
from the sequences of poles diverging to $\iy$.
If $a,b,c,d$ are four reals, or two reals and one pair of complex
conjugates, or two pairs of complex conjugates such that
$|ab|,|ac|,|ad|,|bc|,|bd|,\allowbreak|cd|<1$, then the Askey-Wilson polynomials
are real-valued and their orthogonality can be rewritten as an integral
over $x=\thalf(z+z^{-1})\in[-1,1]$ plus a finite sum over real $x$-values
outside $[-1,1]$. This finite sum does not occur if
$|a|,|b|,|c|,|d|<1$.

Now $k_n$, $h_n$, $B_n$, $D$ and $\la_n$ in \S\ref{sec:43} can be specified
for the Askey-Wilson case as follows.
\begin{equation*}
k_n=2^n(abcdq^{n-1};q)_n,\quad
\frac{h_n}{h_0}=\frac{1-abcdq^{-1}}{1-abcdq^{2n-1}}\,
\frac{(q,ab,ac,ad,bc,bd,cd;q)_n}{(abcdq^{-1};q)_n}\,,
\end{equation*}
\begin{multline}
B_n=
\Bigl((a+b+c+d) (q-a b c d q^{n-1}-a b c d q^n+a b c d q^{2n})\\
+(bcd+abd+acd+abc)(1-q^n-q^{n+1}+a b c d q^{2n-1})\Bigr)\\
\times\frac{q^{n-1}}{2(1-a b c d q^{2 n-2})(1-a b c d q^{2 n})}\,,
\label{eq:35}
\end{multline}
\begin{multline}
\mbox{$Dp_n=\la_n p_n$, where}\\
\thalf(1-q^{-1})(Df)[z]=
v(z)f[qz]-\bigl(v(z)+v(z^{-1})\bigr)f[z]+v(z^{-1})f[q^{-1}z],\\
v(z)=\frac{(1-az)(1-bz)(1-cz)(1-dz)}{(1-z^2)(1-qz^2)}\,,\\
\thalf(1-q^{-1})\la_n=(q^{-n}-1)(1-abcdq^{n-1}).
\label{eq:69}
\end{multline}

Define an operator $L$ acting on symmetric Laurent polynomials:
\begin{align}
&(Lf)[z]:=
\Bigl((1-az)(1-bz)(1-cz)(1-dz)\,z^{-2}\,f[qz]
\nonumber\\
&\qquad
-(1-a/z)(1-b/z)(1-c/z)(1-d/z)\,z^2\,f[q^{-1}z]\Bigr)\,(z-z^{-1})^{-1}
\label{eq:33}\\\noalign{\allowbreak}
&=\thalf(1-q^2)\,
\frac{(az,az^{-1},bz,bz^{-1},cz,cz^{-1},dz,dz^{-1};q^2)_\iy}
{(qz^2,qz^{-2};q^2)_\iy}
\nonumber\\
&\qquad\times\frac{\de_{q^2}}{\de_{q^2}x}
\Biggl(\frac{(qz^2,qz^{-2};q^2)_\iy\,f[z]}
{(qaz,qaz^{-1},qbz,qbz^{-1},qcz,qcz^{-1},qdz,qdz^{-1};q^2)_\iy}\Biggr).
\label{eq:34}
\end{align}
Here
\begin{equation}
\frac{\de_{q}}{\de_{q}x}\,g[z]:=
\frac{2(g[q^\half z]-g[q^{-\half}z])}{(q^\half-q^{-\half})(z-z^{-1})}
\end{equation}
is a divided $q$-difference operator (see \cite[\S5]{05}).
It tends to $\frac d{dx} g(x)$ as $q\uparrow1$.

Then $L$ sends symmetric Laurent polynomials of degree $n$ to
symmetric Laurent polynomials
of degree $n+1$,
\begin{equation}
\ga_n=2(abcd\,q^n-q^{-n}),
\end{equation}
and $L$ is skew symmetric with respect to the inner product
\eqref{eq:30}, so \eqref{eq:24} holds.
For the proof of \eqref{eq:24} note that
\begin{align*}
&\oint_C \frac{(1-az)(1-bz)(1-cz)(1-dz)f[qz]\,g[z]}{z^2(z-z^{-1})}\,
w(z)\,\frac{dz}z\\
\noalign{\allowbreak}
&\;\;=\oint_C \frac{f[qz]\,g[z]}{(qz)^{-2}(qz-(qz)^{-1})}\,
\frac{(q^2z^2,q^{-2}z^{-2};q)_\iy}
{(qaz,az^{-1},qbz,bz^{-1},qcz,cz^{-1},qdz,dz^{-1};q)_\iy}\,\frac{dz}z\\
\noalign{\allowbreak}
&\;\;=\oint_C \frac{f[z]\,g[q^{-1}z]}{z^{-2}(z-z^{-1})}\,
\frac{(z^2,z^{-2};q)_\iy}
{(az,qaz^{-1},bz,qbz^{-1},cz,qcz^{-1},dz,qdz^{-1};q)_\iy}\,\frac{dz}z\\
\noalign{\allowbreak}
&\;\;
=\oint_C \frac{f[z]\,(1-az)(1-bz)(1-cz)(1-dz)g[q^{-1}z]}{z^{-2}(z-z^{-1})}\,
w(z)\,\frac{dz}z\,.
\end{align*}
(Actually, the contour deformation above can be done with avoidance of poles
in the generic case of complex $a,b,c,d$ such that the four line segments
connecting $a,b,c,d$ with 0 avoid the four halflines
$\{ta^{-1},tb^{-1},tc^{-1},td^{-1}\mid t\ge1$.)
Alternatively, we can observe that $L=[D,X]$, $\ga_n=\la_{n+1}-\la_n$
with $D,\la_n$ given by \eqref{eq:69}, and apply Proposition
\ref{th:62}.

By Proposition \ref{th:27} we have the structure relation
\eqref{eq:28}, which can be more explicitly written (with usage of
\eqref{eq:33}) as
\begin{multline}
(Lp_n)[z]=-\frac{(1-abcdq^{n-1})p_{n+1}[z]}{q^n(1-abcdq^{2n-1})}
+(1-abq^{n-1})(1-acq^{n-1})(1-adq^{n-1})
\\
\times(1-bcq^{n-1})
(1-bdq^{n-1})(1-cdq^{n-1})\,
\frac{(1-q^n)p_{n-1}[z]}{q^{n-1}(1-abcdq^{2n-1})}
\label{eq:18}
\end{multline}
We can also write the lowering and raising relations
\eqref{eq:31}, \eqref{eq:32} more explicitly (with usage of
\eqref{eq:33} and \eqref{eq:35}):
\begin{multline}
(Lp_n)[z]-(abcdq^n-q^{-n})(z+z^{-1}-2B_n)\,p_n[z]
=(1-abq^{n-1})
(1-acq^{n-1})
\\
\times(1-adq^{n-1})(1-bcq^{n-1})
(1-bdq^{n-1})(1-cdq^{n-1})\,
\frac{(1+q)(1-q^n)p_{n-1}[z]}{q^n(1-abcdq^{2n-2})}\,,
\label{eq:76}
\end{multline}
\begin{multline}
(Lp_n)[z]+(abcdq^{n-1}-q^{1-n})
(z+z^{-1}-2B_n)\,p_n[z]\\
=-\frac{(1+q)(1-abcdq^{n-1})}{q^n(1-abcdq^{2n})}\,p_{n+1}[z].
\label{eq:77}
\end{multline}

Lowering and raising relations for Askey-Wilson polynomials
were earlier obtained by Bangerezako \cite[(41)]{20}. He obtained his
lowering and raising operator
by the factorization method, see Proposition 1 in \cite{20}.
His lowering and raising relation can be obtained from \eqref{eq:76} and
\eqref{eq:77}, respectively, by adding
$\thalf(1-q^{-1})(z-qz^{-1})(Dp_n)[z]$ to the left-hand sides of these
relations (cf.~Remark \ref{th:75}).
\begin{remark}
It turns out that
the operator \eqref{eq:33}, which occurs in the structure relation for
Askey-Wilson polynomials, is the trigonometric case $p=0$ of the
difference operator $\De(a,b,c,d)$ with elliptic coefficients given by
Rosengren \cite[\S6]{21} (also replace $q$ by $q^2$ in Rosengren's
operator in order to arrive at \eqref{eq:33}).
Rosengren observes, following Rains \cite{23},
that the operators $\De(a,b,c,d)$ with $abcd=q^{-N}$ ($N=0,1,2,\ldots$)
generate a representation of the Sklyanin algebra as in
Sklyanin \cite[Theorem 4]{22}. Indeed, if we write
$L=L_{a,b,c,d}$ for the operator in \eqref{eq:33} then
\begin{equation}
L_{a,b,ce,de^{-1}}L_{qa,qb,q^{-1}c,q^{-1}d}=
L_{a,b,c,d}L_{qa,qb,q^{-1}ce,q^{-1}de^{-1}}.
\end{equation}
This quasi-commutation relation, which can be proved in a straightforward
way, is a specialization ($p=0$, $n=1$) of Rains' relation
\cite[(3.27)]{23} and it is a way to implement the relations in the
Sklyanin algebra.
\end{remark}
\section{Continuous $q$-Jacobi polynomials}
\label{sec:45}
{\em Continuous $q$-Jacobi polynomials}
\[
P_n[z]=P_n^{(\al,\be)}\bigl(\thalf(z+z^{-1})\mid q\bigr)
\]
(see \cite[\S7.5]{06}, \cite[\S3.10]{04})
can be obtained as restrictions of Askey-Wilson polynomials in
two different ways (see \cite[\S4.10]{04}):
\begin{align}
P_n^{(\al,\be)}\bigl(x\mid q\bigr)&=
\frac{q^{(\half\al+\frac14)n}\,
p_n(x;q^{\half\al+\frac14},-q^{\half\be+\frac14},q^{\frac14},-q^{\frac14}\mid
q^\half)}
{(-q^{\half(\al+\be+1)};q^\half)_n\,(q;q)_n}
\label{eq:49}
\\
&=\frac{q^{(\half\al+\frac14)n}\,
p_n(x;q^{\half\al+\frac14},q^{\half\al+\frac34},
-q^{\half\be+\frac14},-q^{\half\be+\frac34}\mid q)}
{(-q^{\half(\al+\be+1)};q^\half)_{2n}\,(q;q)_n}\,.
\label{eq:09}
\end{align}
These polynomials are orthogonal with respect to the inner product
\begin{equation}
\lan f,g\ran:=
\oint f[z]\,g[z]\,
\frac{(z^2,z^{-2};q)_\iy}
{(q^{\half\al+\frac14}z,q^{\half\al+\frac14}z^{-1},
-q^{\half\be+\frac14}z,-q^{\half\be+\frac14}z^{-1};q^\half)_\iy}\,\frac{dz}z\,.
\end{equation}
The coefficients $A_n$ and $C_n$ in \eqref{eq:58} here become:
\begin{align}
A_n&=\frac{(1-q^{n+1})(1-q^{n+\al+\be+1})}
{2 q^{\half\al+\frac14}(1-q^{n+\half(\al+\be+1)})(1-q^{n+\half(\al+\be+2)})}\,,
\\
C_n&=\frac{q^{\half\al+\frac14}(1-q^{n+\al})(1-q^{n+\be})}
{2(1-q^{n+\half(\al+\be)})(1-q^{n+\half(\al+\be+1)})}\,.
\end{align}

Corresponding to \eqref{eq:49}, \eqref{eq:09}
we can obtain two versions of the operators $L$ and $D$
for continuous $q$-Jacobi polynomials by specialization of
\eqref{eq:33} and \eqref{eq:69}. For $L$ this becomes:
\begin{multline}
(Lf)[z]=\frac{v(z) f[q^\half z]-v(z^{-1}) f[q^{-\half}z]}{z-z^{-1}}\,,
\\
v(z)=(1-q^{\half\al+\frac14}z)(1+q^{\half\be+\frac14}z)(1-q^\half z^2)\,
z^{-2},
\end{multline}
\begin{multline}
(\tilde Lf)[z]=\frac{\tilde v(z) f[qz]-\tilde v(z^{-1}) f[q^{-1}z]}
{z-z^{-1}}\,,
\\
\tilde v(z)=(1-q^{\half\al+\frac14}z)(1-q^{\half\al+\frac34}z)
(1+q^{\half\be+\frac14}z)(1+q^{\half\be+\frac34}z)\,z^{-2}.
\end{multline}
Related to these operators are coefficients $\ga_n$ respectively
$\tilde\ga_n$ (see \eqref{eq:57}):
\begin{equation}
\ga_n=2(q^{\half(n+\al+\be+2)}-q^{-\half n}),\qquad
\tilde\ga_n=2(q^{n+\al+\be+2}-q^{-n}).
\end{equation}

The form \eqref{eq:34} of the operator $L$
can also be specialized corresponding to
\eqref{eq:49} and \eqref{eq:09}.
In particular, corresponding to \eqref{eq:49} we obtain:
\begin{multline}
(Lf)[z]=\thalf(1-q)
\frac{(q^{\half\al+\frac14}z,q^{\half\al+\frac14}/z,
-q^{\half\be+\frac14}z,-q^{\half\be+\frac14}/z;q)_\iy}
{(q^{\frac32}z^2,q^{\frac32}z^{-2};q^2)_\iy}
\\
\times\frac{\de_q}{\de_q x}
\left(\frac{(q^\half z^2,q^\half/z^2;q^2)_\iy\,f[z]}
{(q^{\half\al+\frac34}z,q^{\half\al+\frac34}/z,
-q^{\half\be+\frac34}z,-q^{\half\be+\frac34}/z;q)_\iy}\right).
\label{eq:60}
\end{multline}

By \eqref{eq:28} and the formulas above we get two structure relations
for continuous $q$-Jacobi polynomials:
\begin{align}
\frac{v(z) P_n[q^\half z]-v(z^{-1}) P_n[q^{-\half}z]}{z-z^{-1}}
&=\ga_n A_n P_{n+1}[z]-
\ga_{n-1} C_n P_{n-1}[z].
\label{eq:59}
\\
\frac{\tilde v(z) P_n[q z]-\tilde v(z^{-1}) P_n[q^{-1}z]}{z-z^{-1}}
&=\tilde\ga_n A_n P_{n+1}[z]-
\tilde\ga_{n-1} C_n P_{n-1}[z].
\end{align}
As in \eqref{eq:31}, \eqref{eq:32}, or more generally in
Remark \ref{th:75}, lowering and raising relations can also be written
explicitly here (I will skip this).

Ismail \cite[Theorem 15.5.2]{24} gives
an explicit lowering relation for continuous $q$-Jacobi polynomials
which is proved by using a determinant formula due to Christoffel.
This expresses orthogonal polynomials with respect to an orthogonality
measure $\Phi(x)d\mu(x)$ in terms of polyomials with orthogonality measure
$d\mu(x)$ ($\Phi(x)$ a nonnegative polynomial on supp($\mu$)),
see \cite[Theorem 2.7.1]{24}.

The structure relation \eqref{eq:26} for Jacobi
polynomials is a limit case of both structure relations above.
This can be most easily seen in a formal way by using the expressions
for $L$ and $\tilde L$ involving the divided $q$-difference. Let us
give the details for $L$ as given by \eqref{eq:60}. Rewrite \eqref{eq:59} as
\begin{equation}
\tfrac2{1-q}(L_qP_n)[z;q]=\tfrac2{1-q}\ga_n(q) A_n(q) P_{n+1}[z;q]-
\tfrac2{1-q}\ga_{n-1}(q) C_n P_{n-1}[z;q],
\label{eq:61}
\end{equation}
emphasizing the $q$-dependence.
As $q\uparrow1$, $P_n[z;q]\to P_n[z]=P_n^{(\al,\be)}
\bigl(\tfrac{z+z^{-1}}2\bigr)$ (see \cite[(5.10.1)]{04}).
Also, $A_n(q)$, $C_n(q)$ and $\tfrac2{1-q}\ga_n(q)$ tend to $A_n$, $C_n$ and
$4\ga_n$
for Jacobi polynomials as given in \S\ref{sec:44}. Thus the \RHS\
of \eqref{eq:61} tends to the \RHS\ of \eqref{eq:28} for
Jacobi polynomials, i.e., to the \RHS\ of \eqref{eq:26}.
Now consider the \LHS\ of \eqref{eq:61} with \eqref{eq:60} substituted.
Use that
\[
\lim_{q\uparrow1}\frac{(q^a z;q)_\iy}{(z;q)_\iy}=(1-z)^{-a}
\]
(see \cite[(1.3.19)]{06}). Hence
\begin{align*}
\frac{(q^{\half\al+\frac14}z,\tfrac{q^{\half\al+\frac14}}z,
-q^{\half\be+\frac14}z,-\tfrac{q^{\half\be+\frac14}}z;q)_\iy}
{(q^{\frac34}z,\tfrac{q^{\frac34}}z,-q^{\frac34}z,
-\tfrac{q^{\frac34}}z;q)_\iy}&\to
(2-z-\tfrac1z)^{\half-\half\al}(2+z+\tfrac1z)^{\half-\half\be},
\\
\frac{(q^{\frac14}z,\tfrac{q^\frac14}z,-q^{\frac14}z,
-\tfrac{q^{\frac14}}z;q)_\iy}
{(q^{\half\al+\frac34}z,\tfrac{q^{\half\al+\frac34}}z,
-q^{\half\be+\frac34}z,-\tfrac{q^{\half\be+\frac34}}z;q)_\iy}&\to
(2-z-\tfrac1z)^{\half+\half\al}(2+z+\tfrac1z)^{\half+\half\be}.
\end{align*}
Hence
$\tfrac2{1-q}(L_qP_n)[z;q]\to 4(LP_n)[z]$ with $L$ given by
\eqref{eq:38}. So \eqref{eq:61} tends to 4 times \eqref{eq:26}
as $q\uparrow1$.
\section{Continuous $q$-ultraspherical polynomials}
\label{sec:47}
In view of \cite[(7.5.34)]{06} and \eqref{eq:49}
{\em continuous $q$-ultraspherical polynomials} are specializations
of Askey-Wilson polynomials as follows:
\begin{multline}
C_n[z]=C_n\bigl(\thalf(z+z^{-1});t\mid q\bigr)\\
=\frac{(t;q^\half)_n}{(q^\half t;q)_n\,(q;q)_n}\,
p_n\bigl(\thalf(z+z^{-1});
t^\half,-t^\half,q^{\frac14},-q^{\frac14}\mid q^\half\bigr).
\label{eq:50}
\end{multline}
By specialization of \eqref{eq:33} and \eqref{eq:34}
we obtain the two ways of writing the operator $L$ corresponding to
\eqref{eq:50}:
\begin{align}
(Lf)[z]&=\frac{(1-tz^2)(z^{-2}-q^\half)}{z-z^{-1}}\,f[q^\half z]-
\frac{(1-tz^{-2})(z^2-q^\half)}{z-z^{-1}}\,f[q^{-\half}z]\\
&=\thalf(1-q)\,\frac{(tz^2,tz^{-2};q^2)_\iy}
{(q^{\frac32}z^2,q^{\frac32}z^{-2};q^2)_\iy}\,
\frac{\de_q}{\de_q x}\left(
\frac{(q^\half z^2,q^\half z^{-2};q^2)_\iy}
{(qtz^2,qtz^{-2};q^2)_\iy}\,f[z]\right).
\end{align}
Then the structure relation \eqref{eq:18} for Askey-Wilson polynomials
specializes to the following structure relation for the polynomials $C_n$:
\begin{multline}
\frac{(1-tz^2)(z^{-2}-q^\half)}{z-z^{-1}}\,C_n[q^\half z]-
\frac{(1-tz^{-2})(z^2-q^\half)}{z-z^{-1}}\,C_n[q^{-\half}z]\\
=-\,\frac{(1-tq^{n+\half})(1-q^{n+1})}{q^{\half n}(1-tq^n)}\,C_{n+1}[z]
+\frac{(1-tq^{n-\half})(1-t^2q^{n-1})}{q^{\half n-\half}(1-tq^n)}\,C_{n-1}[z].
\label{eq:54}
\end{multline}

In \cite[\S2]{11} some lowering and raising relations for
continuous $q$-ultra\-spherical polynomials were given which were obtained
as $A_1$ cases of the lowering and raising relations for
the $A_n$ Macdonald polynomials in \cite{12}. For instance:
\begin{align}
&-\,\frac{t-z^{-2}}{z-z^{-1}}\,C_n[q^\half z]+
\frac{t-z^2}{z-z^{-1}}\,C_n[q^{-\half}z]+
q^{-\half n}(z+z^{-1})\,C_n[z]
\nonumber\\
&\qquad\qquad\qquad\qquad\qquad\qquad\qquad\qquad
=(q^{-\half n}-t^2q^{\half n-1})\,C_{n-1}[z],
\label{eq:51}
\\
&-\,\frac{tz^2-1}{z-z^{-1}}\,C_n[q^\half z]+
\frac{tz^{-2}-1}{z-z^{-1}}\,C_n[q^{-\half} z]+
q^{-\half n}(z+z^{-1})\,C_n[z]
\nonumber\\
&\qquad\qquad\qquad\qquad\qquad\qquad\qquad\qquad
=(q^{-\half n}-q^{\half n+1})\,C_{n+1}[z].
\label{eq:52}
\end{align}
By subtracting \eqref{eq:51} from \eqref{eq:52} we obtain another
structure relation for continuous $q$-ultraspherical polynomials
of the form \eqref{eq:23}, but with an operator $L$ which is not
skew symmetric:
\begin{multline}
z^{-1}(1-tz^2)C_n[q^\half z]+z(1-tz^{-2})C_n[q^{-\half}z]\\
=q^{-\half n}(1-q^{n+1})C_{n+1}[z]-q^{-\half n}(1-t^2 q^{n-1})C_{n-1}[z].
\label{eq:53}
\end{multline}

Relation \eqref{eq:53} is connected with relation
\eqref{eq:54} by yet another
relation of the form \eqref{eq:23}:
\begin{multline}
-\,\frac{(1-tz^2)((1+z^{-2})}{z-z^{-1}}\,C_n[q^\half z]
+\frac{(1-tz^{-2})(1+z^2)}{z-z^{-1}}\,C_n[q^{-\half}z]\\
=
\frac{(q^{-\half n}+tq^{\half n})(1-q^{n+1})}{1-tq^n}\,C_{n+1}[z]+
\frac{(q^{-\half n}+tq^{\half n})(1-t^2q^{n-1})}{1-tq^n}\,C_{n-1}[z].
\label{eq:55}
\end{multline}
Indeed,
\eqref{eq:54} equals $\thalf(q-1)$ times \eqref{eq:55} plus
$\thalf(q+1)$ times \eqref{eq:53}.

Relation \eqref{eq:55} is a consequence of
the three-term recurrence relation
\begin{equation}
(z+z^{-1})C_n[z]=\frac{1-q^{n+1}}{1-tq^n}\,C_{n+1}[z]
+\frac{1-t^2q^{n-1}}{1-tq^n}\,C_{n-1}[z]
\end{equation}
(see \cite[(3.10.17)]{04}) and the second order $q$-difference formula
\begin{equation}
\frac{1-tz^2}{1-z^2}\,C_n[q^\half z]+
\frac{1-tz^{-2}}{1-z^{-2}}\,C_n[q^{-\half}z]=
(q^{-\half n}+q^{\half n} t)\,C_n[z],
\end{equation}
which last formula is the specialization by \eqref{eq:50} of the
second order $q$-difference formula \eqref{eq:69} for
Askey-Wilson polynomials.
\section{Big $q$-Jacobi polynomials}
\label{sec:46}
{\em Big $q$-Jacobi polynomials}
(see \cite[\S7.3]{06}, \cite[\S3.5]{04})
are given by
\begin{equation}
P_n(x)=P_n(x;a,b,-c;q):=
\qhyp32{q^{-n},abq^{n+1},x}{aq,-cq}{q,q}
\end{equation}
They can be obtained as limits of Askey-Wilson polynomials
(see \cite[(4.1.3)]{04}):
\begin{multline}
\label{eq:17}
\lim_{\ep\to0}\frac{\ep^n}{(aq,-cq,-\ep^2 bc^{-1};q)_n}\,
p_n\bigl(\thalf(a^{-1}x+ax^{-1});\ep,\ep^{-1}aq,-\ep^{-1}cq,-\ep bc^{-1}
\mid q\bigr)\\
=P_n(x;a,b,-c;q).
\end{multline}
Make in \eqref{eq:18} together with \eqref{eq:19} the substitutions
corresponding to the \LHS\ of \eqref{eq:17} and let $\ep\to0$. Then we
we obtain a structure relation for big $q$-Jacobi polynomials:
\begin{multline}
(LP_n)(x)=
\frac{(1-aq^{n+1})(1+cq^{n+1})(1-abq^{n+1})}{q^{n+2}ac(1-abq^{2n+1})}\,
P_{n+1}(x)\\
-\frac{(1-q^n)(1-bq^n)(1+abc^{-1}q^n)}{1-abq^{2n+1}}\,
P_{n-1}(x),
\label{eq:40}
\end{multline}
where
\begin{equation}
(Lf)(x):=\frac{(1-x)(1+bc^{-1}x)\,f(qx)
-(1-a^{-1}q^{-1}x)(1+c^{-1}q^{-1}x)\,f(q^{-1}x)}x.
\label{eq:56}
\end{equation}
Note that we can rewrite \eqref{eq:56} as
\begin{equation}
(Lf)(x)=(q-q^{-1})\,\frac{(x,-bc^{-1}x;q^2)_\iy}
{(qa^{-1}x,-qc^{-1}x;q^2)_\iy}\,
d_{q,x}\left(
\frac{(a^{-1}x,-c^{-1}x;q^2)_\iy}{(qx,-qbc^{-1}x;q^2)_\iy}\,f[x]
\right),
\end{equation}
where $d_{q,x}$ is a central $q$-derivative:
\begin{equation}
d_{q,x}\bigl(g(x)\bigr):=\frac{g(qx)-g(q^{-1}x)}{(q-q^{-1})x}\,.
\end{equation}

We can rewrite \eqref{eq:40} in two steps as the explicit
structure relation of the form
\begin{equation}
(x-1)(bx+c)\,D_{q,x}\bigl(P_n(x)\bigr)=
\tilde a_n P_{n+1}(x)+\tilde b_n P_n(x)+\tilde c_n P_{n-1}(x)
\label{eq:41}
\end{equation}
given in \cite{09}, where the $q$-derivative $D_q$ is defined by
\eqref{eq:39}.
First eliminate $P_n(q^{-1}x)$ from \eqref{eq:40}
by means of the second order $q$-difference equation for big
$q$-Jacobi polynomials. Then we obtain from \eqref{eq:40} a relation of
the form
\begin{equation}
(x-1)(bx+c)\,D_{q,x}\bigl(P_n(x)\bigr)=
\al_n P_{n+1}(x)+(\de_n x+\be_n)P_n(x)+\ga_n P_{n-1}(x).
\label{eq:42}
\end{equation}
Next eliminate $x P_n(x)$ from \eqref{eq:42} by means of the
three-term recurrence relation in order to obtain a relation of
the form \eqref{eq:41}.
%
%:References

\end{document}